\documentclass[a4apper,10pt]{amsart}

\usepackage{latexsym}
\usepackage{amsmath}
\usepackage{amsthm}
\usepackage{amssymb}
\usepackage{amsfonts}
\usepackage{amsxtra}
\usepackage{setspace}

\newtheorem{thm}{Theorem}
\newtheorem{prop}[thm]{Proposition}
\newtheorem{cor}[thm]{Corollary}
\newtheorem{lem}[thm]{Lemma}

\theoremstyle{definition}
\newtheorem{defn}[thm]{Definition}

\numberwithin{thm}{section}

\newcommand{\Z}{\mathbb{Z}}
\newcommand{\N}{\mathbb{N}}
\newcommand{\Q}{\mathbb{Q}}

\newcommand{\C}{\mathbb{C}}

\renewcommand{\Im}{\textrm{Im}}
\renewcommand{\ker}{\textrm{Ker}}
\newcommand{\Ker}{\textrm{Ker}}

\newcommand{\End}{\textrm{End}}
\newcommand{\Verma}{\Delta_{\mathbf{c}}(E;p,q)}
\newcommand{\Irr}{\mathsf{Irr}}

\begin{document}
\title{On singular Calogero-Moser Spaces}
\author{Gwyn Bellamy}
\address{School of Mathematics and Maxwell Institute for Mathematical Sciences, University of Edinburgh, James Clerk Maxwell Building, Kings Buildings, Mayfield Road, Edinburgh EH9 3JZ, Scotland}
\email{G.E.Bellamy@sms.ed.ac.uk}
\maketitle
\begin{abstract}
\noindent Using combinatorial properties of complex reflection groups we show that if the group $W$ is different from the wreath product $\mathfrak{S}_n\wr \mathbb{Z}/m\mathbb{Z}$ and the binary tetrahedral group (labelled $G(m,1,n)$ and $G_4$ respectively in the Shephard-Todd classification), then the generalised Calogero-Moser space $X_{\mathbf{c}}$ associated to the centre of the rational Cherednik algebra $H_{0,\mathbf{c}}(W)$ is singular for all values of the parameter $\mathbf{c}$. This result and a theorem of Ginzburg and Kaledin imply that there does not exist a symplectic resolution of the singular symplectic variety $\frak{h} \times \frak{h}^* / W$ when $W$ is a complex reflection group different from $\mathfrak{S}_n\wr \mathbb{Z}/m\mathbb{Z}$ and the binary tetrahedral group (where $\frak{h}$ is the reflection representation associated to $W$). Conversely it has been shown by Etingof and Ginzburg that $X_{\bf{c}}$ is smooth for generic values of $\bf{c}$ when $W \cong \mathfrak{S}_n \wr \mathbb{Z}/m\mathbb{Z}$. We show that this is also the case when $W$ is the binary tetrahedral group. A theorem of Namikawa then implies the existence of a symplectic resolution in this case, completing the classification. Finally, we note that the above results together with work of Chlouveraki are consistent with a conjecture of Gordon and Martino on block partitions in the restricted rational Cherednik algebra.
\end{abstract}

\section{Introduction}
Let $W$ be an irreducible complex reflection group and $\frak{h}$ its reflection representation. Etingof and Ginzburg \cite{1} associated to $W$ a family of algebras, the \textit{rational Cherednik algebras} $H_{t,\mathbf{c}}(W)$, depending on parameters $t$ and $\mathbf{c}$. The definition is given in Section 2. When $t = 0$, these algebras have large centres and the geometry of the centre strongly influences the representation theory of the algebra. The affine variety $X_{\mathbf{c}}$ corresponding to the centre of the rational Cherednik algebra was called the generalised Calogero-Moser space at $\mathbf{c}$ by Etingof and Ginzburg. They showed \cite[Corollary 1.14]{1}, that for generic values of the parameter $\mathbf{c}$, $X_{\mathbf{c}}$ is smooth when $W \cong G(m,1,n)$. However, Gordon \cite[Proposition 7.3]{6} showed that, for many Weyl groups $W$ not of type $A$ or $B(=C)$, $X_{\mathbf{c}}$ is a singular variety for all choices of the parameter $\mathbf{c}$. Using similar methods we extend this result to all irreducible complex reflection groups. 

\begin{thm}\label{thm:intro}
Let $W$ be an irreducible complex reflection group, not isomorphic to $G(m,1,n)$ or $G_4$, and $X_{\mathbf{c}}$  the generalised Calogero-Moser space associated to $W$. Then $X_{\mathbf{c}}$ is a singular variety for all choices of the parameter $\mathbf{c}$. Conversely for $W \cong G_4$, $X_{\mathbf{c}}$ is a smooth variety for generic values of $\mathbf{c}$.
\end{thm}

\noindent This completes the classification of rational Cherednik algebras for which $X_{\mathbf{c}}$ is smooth for generic $\mathbf{c}$. \\

In \cite[Corollary 1.21]{8}, Ginzburg and Kaledin show that the existence of a symplectic resolution of the symplectic singularity $\frak{h} \times \frak{h}^* / W$ implies that $X_{\mathbf{c}}$ is smooth for generic $\mathbf{c}$. This result, together with Theorem \ref{thm:intro} above implies the following corollary.
\begin{cor}
Let $W$ be an irreducible complex reflection group with reflection representation $\frak{h}$. Then there does not exist a symplectic resolution of $\frak{h} \times \frak{h}^* / W$ when $W \not\cong G(m,1,n)$ or $G_4$.
\end{cor}
It has been show by Wang \cite[Proposition 1]{wang}, that there exists a symplectic resolution of \mbox{$\frak{h} \times \frak{h}^* / W$} when $W \cong G(m,1,n)$, for all $m$ and $n$. Similarly, since $X_{\bf{c}}$ is smooth for generic values of $\bf{c}$ when $W \cong G_4$, a result of Namikawa \cite[Corollary 2.10]{9}, implies 
\begin{cor}
There exists a symplectic resolution of the singular symplectic variety $\frak{h} \times \frak{h}^* / G_4$.
\end{cor}
In order to prove Theorem \ref{thm:intro} we show that the restricted rational Cherednik algebra $\bar{H}_{0,\bf{c}}(W)$ has irreducible representations of dimension $< |W|$ for all values of $\bf{c}$ when $W$ is different from $G(m,1,n)$ and $G_4$. This implies that there exist blocks in $\bar{H}_{0,\bf{c}}(W)$ with nonisomorphic irreducible modules. Therefore the corresponding block partition of $\Irr(W)$, as described in \cite{12}, is trivial for generic values of $\bf{c}$ if and only if $W$ is $G(m,1,n)$ or $G_4$. A conjecture of Gordon and Martino \cite{12} then implies that the partitioning of $\Irr(W)$ induced by the Rouquier families of the Hecke algebra $\mathcal{H}_{\bf{q}}(W)$ should also be trivial for generic choices of $\bf{c}$ if and only if $W$ is $G(m,1,n)$ or $G_4$. Work of Chlouveraki \cite{13} on the cyclotomic Hecke algebras of exceptional complex reflection groups shows that this is indeed the case.

\section{The rational Cherednik algebra at $t = 0$}

\subsection{Definitions and notation}\label{subsection:defns}

Let $W$ be a complex reflection group, $\frak{h}$ its reflection
representation over $\C$ with $\dim ( \frak{h} ) = n$, and $\mathcal{S}$ the set of all complex reflections in $W$. Let $\omega : \frak{h} \oplus \frak{h}^* \rightarrow \C$ be the symplectic form on $\frak{h} \oplus \frak{h}^*$ given by
$\omega((f_1,f_2),(g_1,g_2)) = f_2(g_1) - g_2(f_1)$ and $\mathbf{c} : \mathcal{S} \rightarrow \C$ a
$W$-invariant function. For $s \in S$, define $\omega_s : \frak{h}
\oplus \frak{h}^* \rightarrow \C$ to be the restriction of $\omega$ on
$\Im (1 - s)$ and the zero form on $\Ker (1 - s)$. The \textit{rational Cherednik algebra} at parameter $t = 0$, as introduced by Etingof and Ginzburg \cite[page 250]{1}, is the quotient of the skew group algebra of the tensor algebra $T(\frak{h} \oplus \frak{h}^*)$ with $W$, $T(\frak{h} \oplus \frak{h}^*) \rtimes W$, by the ideal generated by the relations

\begin{equation}\label{eq:rel}
[x,y] = \sum_{s \in S} \mathbf{c}(s) \omega_s(x,y) s \qquad \qquad
\forall x,y \in \frak{h} \oplus \frak{h}^* 
\end{equation}

\noindent Let $Z_{\mathbf{c}}$ denote the centre of $H_{0,\mathbf{c}}$
and $X_{\mathbf{c}} = \textrm{maxspec}(Z_{\mathbf{c}})$ the affine
variety corresponding to $Z_{\mathbf{c}}$. The space $X_{\mathbf{c}}$ is called the \textit{generalised Calogero-Moser space} associated to the complex
reflection group $W$ at parameter $\mathbf{c}$. By \cite[Proposition 4.5]{1}, we have an inclusion $A = \C[\frak{h}]^W\otimes\C[\frak{h}^*]^W \subset Z_{\mathbf{c}}$ and correspondingly a surjective morphism $\Upsilon_{\mathbf{c}} : X_{\mathbf{c}} \rightarrow \frak{h}/W \times \frak{h}^* / W$. This allows us to define the restricted rational Cherednik algebra $\bar{H}_{0,\mathbf{c}}(W)$ as

\begin{displaymath}
\bar{H}_{0,\mathbf{c}}(W) := \frac{H_{0,\mathbf{c}}(W)}{\langle A_+ \rangle}
\end{displaymath}

\noindent where $A_+$ denotes the ideal in $A$ of elements with zero constant term. From the defining relations (\ref{eq:rel}) we see that putting $\frak{h}^*$ in degree $1$, $\frak{h}$ in degree $-1$ and $\C W$ in degree $0$ defines a $\Z$-grading on the rational Cherednik algebra $H_{t,\bf{c}}(W)$. The ideal $\langle A_+ \rangle$ is generated by elements that are homogeneous with respect to this grading, therefore $\bar{H}_{0,\bf{c}}$ is also a $\Z$-graded algebra.\\ 

\noindent Let $\Irr(W)$ be a complete set of non-isomorphic irreducible representation of $W$. We denote by $\C[\frak{h}]^{coW}$ the coinvariant ring $\C[\frak{h}]/ C[\frak{h}]^W_+$, where $\C[\frak{h}]^W_+$ is the ideal in $\C[\frak{h}]$ generated by the elements in $\C[\frak{h}]^W$ with zero constant term. We follow the notation introduced in \cite{6} and define 
\begin{displaymath}
M(\lambda) := \bar{H}_{0,\mathbf{c}} \otimes_{\C[\frak{h}]^{coW} \rtimes W} \lambda 
\end{displaymath}
to be the baby Verma $\bar{H}_{0,\mathbf{c}}$-module associated to $\lambda$. This module is a graded $\bar{H}_{0,\mathbf{c}}$-module with $M(\lambda)_i = 0$ for $i > 0$. By \cite[Proposition 4.3]{6}, $M(\lambda)$ has a simple head which we denote $L(\lambda)$.\\

\noindent We follow the notation of \cite{7} with regards to complex reflection groups, and set $d = m/p$ when considering the groups $G(m,p,n)$. For an arbitrary finite dimensional $\Z$-graded vector space $M = \oplus_{i \in \Z} M_i$, the Poincar\'e polynomial of $M$ will be denoted $P(M,t)$. We denote by $f_{\lambda}(t)$ the \textit{fake polynomial} of the $\lambda \in \Irr(W)$. This is defined as
\begin{displaymath}
f_{\lambda}(t) := \sum_{i \in \Z_{\ge 0}} (\C[\frak{h}]^{coW}_i : \lambda)t^i
\end{displaymath}
where $(\C[\frak{h}]^{coW}_i : \lambda)$ is the multiplicity of $\lambda$ in $i^{th}$ degree of the coinvariant ring $\C[\frak{h}]^{coW}$ (thought of here as a graded $W$-module).\\

\noindent We will also require the surjective map $\Theta : \Irr(W) \rightarrow \Upsilon^{-1}(0)$, taking $\lambda$ to the annihilator of $L(\lambda)$ in $Z_{\mathbf{c}}$, as defined in \cite[paragraph 5.4]{6}. This map has the property that a fiber $\Theta^{-1}(\frak{m})$ is a singleton set if and only if $\frak{m}$ is a smooth closed point in $X_{\bf c}$ (\cite[Theorem 5.6]{6}).

\subsection{General results}
Let $\{ s_1, \dots s_k \}$ be a conjugacy class consisting of complex reflections in $W$ and $\zeta$ the eigenvalue of $s_1$ (and hence all $s_i$) not equal to $1$ when thinking of $W$ as a subgroup of $GL(\frak{h})$. For $1 \le i\le k$, let $\omega_{s_i}$ be the restricted symplectic form on $\frak{h} \oplus \frak{h}^*$ as defined above. Let $\pi_{s_i} : \frak{h} \oplus \frak{h}^* \rightarrow \Im ( 1 - s_i )$ be the projection map along $\ker (1 - s_i)$, so that $\omega_{s_i} = \omega \circ \pi_{s_i}$, and define $\Omega = \sum_{i = 1}^k \omega_{s_i}$. 

\begin{lem}\label{lem:nondeg}
Let $W$, $\omega$ and $\Omega$ be as above. Then 
\begin{displaymath}
\Omega = \frac{k}{n}(1 - \zeta)^{-1}(1 - \zeta^{-1})^{-1}(2 - \zeta - \zeta^{-1}) \omega.
\end{displaymath}
\end{lem}

\begin{proof}
Since each $\omega_{s_i}$ is alternating and $\C$-linear, $\Omega \in
\bigwedge^2(\frak{h} \oplus \frak{h}^*)$. Let $x \in \frak{h} \oplus
\frak{h}^*$. Then $x$ decomposes uniquely as $x_1 + x_2$, with $x_1 \in \Im ( 1 - s_i)$ and $x_2 \in \ker (1 - s_i)$. By definition, there exists $y \in \frak{h} \oplus \frak{h}^*$
such that $(1 - s_i)y = x_1$. Then $(1 - gs_ig^{-1})(gy) = g(1 -
s_i)g^{-1}gy = g(1 - s_i)y = gx_1$ implying that $gx_1 \in \Im( 1 -
gs_ig^{-1})$. Also $(1 - s_i)x_2 = 0$ implies that $(1 -
gs_ig^{-1})gx_2 = 0$ hence $gx$ decomposes as $gx_1 + gx_2$ with $gx_1
\in \Im(1 - gs_ig^{-1})$ and $gx_2 \in \ker(1 - gs_ig^{-1})$. Therefore
$\pi_{gs_ig^{-1}} = g \pi_{s_i} g^{-1}$ and $\omega_{s_i} (g^{-1} x, g^{-1} y) = \omega_{gs_ig^{-1}}(x,y)$. Hence $\Omega \in (\bigwedge^2 (\frak{h}^*
\oplus \frak{h}))^W$. By \cite[Lemma 2.23]{1} $\dim (\bigwedge^2 (\frak{h}^*
\oplus \frak{h}))^W = 1$, therefore there exists $\lambda \in \C$ such
that $\Omega = \lambda \omega$. Consider $\Omega'(x,y) = \Omega((x,0),(0,y))$, where $x \in \frak{h}$
and $y \in \frak{h}^*$. Recall that $\zeta$ is the eigenvalue of $s_i$ not equal to $1$, then $\pi_{s_i}(x) = (1 - \zeta)^{-1}(1 - s_i)x$ and $\pi_{s_i}(y) = (1 - \zeta^{-1})^{-1}(1 - s_i)y$. Expanding
$\Omega'(x,y)$
\begin{displaymath}
\Omega'(x,y) = \sum_{i = 1}^k \omega((1 - \zeta)^{-1}(1 - s_i)x,(1 -
\zeta^{-1})^{-1}(1 - s_i)y)
\end{displaymath}
\begin{displaymath}
= (1 - \zeta)^{-1}(1 - \zeta^{-1})^{-1}
\sum_{i = 1}^k \left[\omega(x,y) - \omega(s_ix,y) - \omega(x,s_iy) +
\omega(s_ix,s_iy)\right]
\end{displaymath}
\begin{displaymath}
= (1 - \zeta)^{-1}(1 - \zeta^{-1})^{-1} \omega(x, (\sum_{i = 1}^k 2 - s_i - s_i^{-1})y)
\end{displaymath} 

\noindent Define $\phi = (\sum_{i = 1}^k 2 - s_i - s_i^{-1}) : \frak{h}^*
\rightarrow \frak{h}^*$, a $W$-homomorphism. The trace of $\phi$
is $2nk - (n-1)k - k\zeta - (n-1)k - k\zeta^{-1} = k(2 - \zeta -
\zeta^{-1})$. Since $\frak{h}^*$ is irreducible, Schur's lemma says
that $\phi(y) = \frac{k}{n}(2 - \zeta - \zeta^{-1})y$ and hence $\lambda = \frac{k}{n}(1 - \zeta)^{-1}(1 - \zeta^{-1})^{-1}(2 - \zeta - \zeta^{-1})$. 
\end{proof}

\noindent We also require the notion of a generalised baby Verma module, which are baby Verma modules above points other than the origin in $\frak{h}/W \times \frak{h}^*/W$.

\begin{defn}
Let $(p,q) \in \frak{h}/W \times \frak{h}^*$, $W_q$ the
stabiliser subgroup of $q$ in $W$ and $E$ an irreducible $W_q$-module. Then we define the \textit{generalised baby Verma} module 
\begin{displaymath}
\Verma := H_{0,\mathbf{c}}(W) \otimes_{\C[\frak{h}]^W \otimes
  \C[\frak{h}^*] \rtimes W_q} E
\end{displaymath}
where the action of $\C[\frak{h}]^W \otimes \C[\frak{h}^*] \rtimes W_q$ on $E$ is given by $(f \otimes g \otimes w) \cdot e = f(p)g(q)w \cdot e$ for all $f \in \C[\frak{h}]^W$, $g \in \C[\frak{h}^*]$, $w \in W_q$.
\end{defn}

\noindent Since $\C[\frak{h}]^W \otimes \C[\frak{h}^*]^W \subseteq
Z_{\mathbf{c}}$, Schur's lemma implies that, for every irreducible
$H_{0,\mathbf{c}}$-module $L$, there exists $(p,r ) \in \frak{h}/W \times \frak{h}^*/W$ such that $(f \otimes g) \cdot l = f(p)g(r)l$, for all $l \in
L$, $f,g \in \C[\frak{h}]^W \otimes \C[\frak{h}^*]^W$. Choosing a
point $q$ in the orbit represented by $r$ we write $(p,r) = (p,Wq)$ and say that the irreducible $H_{0,\mathbf{c}}$-module $L$ \textit{lies above} $(p,Wq)$.

\begin{lem}\label{lem:Verma}
Let $L$ be an irreducible $H_{0,\mathbf{c}}$-module lying above $(p,Wq)$. Then there exist $E \in \Irr(W_q)$ and a
surjective $H_{0,\mathbf{c}}$-homomorphism $\phi : \Verma \twoheadrightarrow L$.
\end{lem}

\begin{proof}
The action on $L$ of the commutative ring $\C[\frak{h}^*]$ gives a
decomoposition $L = \oplus_{q' \in \frak{h}^*} L_{q'}^{gen}$ of $L$
into generalised eigenspaces. That is, for each $l \in L_{q'}^{gen}$
and $f \in \C[\frak{h}^*]$, there exists an $N \in \N$ such that $(f -
f(q'))^N \cdot l = 0$ (since $L$ is finite dimensional, we can choose $N$ to be independent of $f$ and $l$).

Choose $q'$ such that $L_{q'}^{gen} \neq 0$, so that $(f - f(q'))^N$ acts as zero on $L_{q'}^{gen}$ for all $f \in \C[\frak{h}^*]^W$. As $L$ lies over $(p,Wq)$ we see that $(f - f(q))$ also acts nilpotently on $L_{q'}^{gen}$ and $f(q) = f(q')$. Since $W$ is a finite group, each orbit in $\frak{h}^*$ is closed, therefore $q' \in Wq$ and we can find $w \in W$ such that $w \cdot q = q'$. Now let $0 \neq L_{q'} \subseteq  L_{q'}^{gen}$ be the space of elements $l$ in $L_{q'}^{gen}$ such that $(f - f(q')) \cdot l = 0$, for all $f \in \C[\frak{h}^*]$. Then $w^{-1}(L_{q'}) \neq 0$ and  $f \cdot (w^{-1}l) =  w^{-1} \cdot ({}^{w}f)(q')l  = f(q)w^{-1} \cdot l$ implies that $w^{-1}(L_{q'}) \subseteq L_q$. Thus $L_q$ is a nonzero $W_q$-submodule of $L$ because $f \cdot (v \cdot l) = v \cdot f(q)l = f(q)(v \cdot l)$ for all $f \in \C[\frak{h}]$, $v \in W_q$ and $l \in L_q$. Choose an irreducible $W_q$-submodule $E$ of $L_q$. The inclusion $E \hookrightarrow L$ induces a $H_{0,\mathbf{c}}$-homomorphism $\phi : \Verma \rightarrow L$. The fact that $L$ is irreducible implies that this is a
surjection.
\end{proof}

\section{Singular generalised Calogero-Moser Spaces}\label{sec:main}

\subsection{The main result}

\begin{thm}\label{thm:main}
For all $W$ not isomorphic to $G(m,1,n)$ or $G_4$ and for all parameters $\mathbf{c}$, the variety $X_{\mathbf{c}}$ is singular.
\end{thm}

\noindent By \cite[Proposition 3.8]{1} the statement of Theorem \ref{thm:main} is equivalent to the statement: \textit{for $W$ not isomorphic to $G(m,1,n)$ or $G_4$ and for all parameters $\mathbf{c}$ there exists an irreducible $H_{0,\bf{c}}(W)$-module $L$ with $\dim L < |W|$}. Therefore Theorem \ref{thm:main} follows from

\begin{prop}\label{prop:mainprop}
For each $W$ not isomorphic to $G(m,1,n)$ or $G_4$, there exists an irreducible $W$-module $\lambda$ such that for all parameters $\bf{c}$, the irreducible $\bar{H}_{0,\bf{c}}(W)$-module $L(\lambda)$ has dimension $< |W|$.
\end{prop}

\noindent The proof of Proposition \ref{prop:mainprop} will occupy the remainder of Section \ref{sec:main}. The irreducible complex reflection groups were classified by Shephard and Todd \cite{5} and either belong to an infinite family labelled $G(m,p,n)$, where $m,p,n \in \N$ and $p \, | \, m$, or to one of 34 exceptional groups $G_4, \dots,  G_{37}$. 

\begin{lem}\label{lem:poly}
Let $W$ be a complex reflection group. Let $\lambda \in \Irr(W)$ be the
unique representation corresponding to a smooth point of $\Upsilon^{-1}(0)$ in $X_{\mathbf{c}}$ i.e. $\Theta(\lambda)$ is smooth in $X_{\mathbf{c}}$. Then the Poincar\'e polynomial of $L(\lambda)$ as a graded vector space is given by

\begin{equation}\label{eq:poin}
P(L(\lambda),t) = \frac{\dim (\lambda)t^{b_{\lambda^*}}P(\C[\frak{h}]^{coW},t)}{f_{\lambda^*}(t)}
\end{equation}

\noindent where $\lambda^*$ is the irreducible $W$-module dual to $\lambda$, and $b_{\lambda}$ the trailing degree of the fake polynomial $f_{\lambda}(t)$.
\end{lem}

\begin{proof}
By \cite[Lemma 4.4, paragraphs (5.2) and (5.4)]{6}, the graded composition factors of $M(\lambda)$ are all of the form $L(\lambda)[i]$, for some $i \ge 0$. Therefore we can find a multiset $\{ i_1, \dots i_k \}$ such that as a graded $W$-module
\begin{displaymath}
M(\lambda) \cong L(\lambda)[i_1] \oplus L(\lambda)[i_2] \oplus \dots \oplus L(\lambda)[i_k].
\end{displaymath}
Since $\Theta(\lambda)$ is a smooth point in $X_{\mathbf{c}}$, \cite[Theorem 1.7]{1} says that $L(\lambda) \cong \C W$ as a $W$-module. Hence it contains a unique copy of the trivial representation $T$. Assume this copy occurs in degree $a$ in
$L(\lambda)$. Then it will occur in degree $a - i_j$ in
$L(\lambda)[i_j]$. As a graded $W$-module, $M(\lambda) \cong
\C[\frak{h}]^{coW}\otimes \lambda$. The fact that $[\tau \otimes
  \lambda : T] = \delta_{\tau\lambda^*}$ implies that the graded
multiplicity of $T$ in $M(\lambda)$ equals the graded multiplicity of
$\lambda^*$ in $\C[\frak{h}]^{coW}$. The graded multiplicity of
$\lambda^*$ in $\C[\frak{h}]^{coW}$ is $f_{\lambda^*}(t)$. Hence
$P(M(\lambda),t) = t^{-a}f_{\lambda^*}(t)P(L(\lambda),t)$. The lowest
nonzero term of $P(L(\lambda),t)$ occurs in degree zero implying that $a = b_{\lambda^*}$. The formula follows by noting that $P(M(\lambda),t)$ is $\dim (\lambda) P(\C[\frak{h}]^{coW},t)$.
\end{proof} 

\noindent Since $L(\lambda)$ is a finite dimensional module, the above lemma shows that the right hand side of equation (\ref{eq:poin}) is a polynomial in $\Z[t,t^{-1}]$ with integer coefficients. Moreover, \cite[Lemma 4.4]{6} shows that it is actually in $\Z[t]$ and that the degree 0 coefficient is $\dim \lambda$.\\

\subsection{The infinite series}
We show that for $p \neq 1$ and $W = G(m,p,n) \neq G(2,2,3)$ we can choose an irreducible representation $\lambda$ of $G(m,p,n)$ such that Lemma \ref{lem:poly} does not hold. Thus $L(\lambda)$ will have dimension $< |G(m,p,n)|$, proving Proposition \ref{prop:mainprop} in this case. The group $G(2,2,3)$ is the Weyl group corresponding to the Dynkin diagram $D_3 = A_3$ and hence $G(2,2,3) \cong S_4$. By \cite[Corollory 16.2]{1}, $X_{\mathbf{c}}$ is smooth for generic and hence all non-zero $\mathbf{c}$ in this case.\\

\noindent We give a description of the parameterization of irreducible $G(m,p,n)$-modules. The reader should consult \cite[pages 379-381]{7} for details. An $m$-multipartition $\underline{\lambda}$ of $n$ is an ordered $m$-tuple of partitions $(\lambda^0, \dots ,\lambda^{m-1})$ such that $|\lambda^0| + \dots + |\lambda^{m-1}| = n$. Let $\mathcal{P}(m)$ denote the set of all $m$-multipartitions of $n$. The cyclic group $\Z / p \Z = \langle g \rangle$ acts on $\mathcal{P}(m)$: $g$ moves each entry of $\underline{\lambda}$ $d$ places to the right i.e.
\begin{displaymath}
g \cdot (\lambda^0, \dots, \lambda^{m-1}) = (\lambda^{m-d},\lambda^{m-d+1}, \dots, \lambda^{m-1},\lambda^0, \dots, \lambda^{m-1}),
\end{displaymath}
(recall from subsection \ref{subsection:defns} that $d = m/p$). For $\underline{\lambda} \in \mathcal{P}(m)$, we denote the orbit $\Z/p\Z \cdot \underline{\lambda}$ by $\{ \underline{\lambda} \}$ and $\textrm{Stab}_{\, \Z/p\Z} \, (\underline{\lambda})$ is the stabliser subgroup of $\Z/p\Z$ with respect to $\underline{\lambda}$. Then the irreducible representations of $G(m,p,n)$ are labelled by distinct pairs $(\{ \underline{\lambda} \},\epsilon)$, where $\epsilon \in \textrm{Stab}_{\Z/p\Z}(\underline{\lambda})$.\\

\noindent Let $(t)_{(n)} = (1 - t) \dots (1 - t^{n-1})(1 - t^n)$ and for $\lambda$ a partition of $n$, denote by $n(\lambda) = \sum_{i = 1}^k (i-1)\lambda_i$ the partition statistic. The young diagram $D_{\lambda}$ of a partition $\lambda$ is the finite subset of $\N \times \N$, justified to the south west (in the French style), representing $\lambda$. For $(i,j) \in D_\lambda$, we denote by $h(i,j)$ the hook length at $(i,j)$. The hook polynomial is defined to be 
\begin{displaymath}
H_{\lambda}(t) = \prod_{(i,j) \in D_\lambda} (1 - t^{h(i,j)}).
\end{displaymath}
\noindent \cite[Corollary 6.4]{7} states that the fake polynomial of the irreducible representation labelled by $(\{ \underline{\lambda} \},\epsilon)$ is

\begin{equation}\label{eq:fake}
f_{\{ \underline{\lambda} \} }(t) = \frac{1 - t^{dn}}{1 - t^{mn}} R_{\{ \underline{\lambda} \} }(t) I_{\underline{\lambda}}(t^m),
\end{equation}

\noindent where 

\begin{displaymath}
R_{\{ \underline{\lambda} \} }(t) = \sum_{\underline{\mu} \in \{ \underline{\lambda} \} } t^{r(\underline{\mu})} \quad \text{ with } \quad r(\underline{\mu}) = \sum_{i = 0}^{m-1} i|\mu^i|, \qquad \text{ and } \quad I_{\underline{\lambda}}(t) = (t)_{(n)} \prod_{i = 1}^m \frac{t^{n(\lambda^i)}}{H_{\lambda^i}(t)}.
\end{displaymath}

\noindent Note that the formula only depends on the orbit and not on the choice of stabiliser.\\

\noindent We wish to calculate the rational function (\ref{eq:poin})
for a well chosen representation $(\{ \underline{\mu} \},\epsilon)$ of the irreducible representations of $G(m,p,n)$. By \cite[Theorem 3.15]{4}, the Poincar\'e polynomial of the coinvariant ring of $W$ is given by

\begin{displaymath}
P(\C[\frak{h}^*]^{coW},t) = \prod_{i = 1}^{n} \frac{1 - t^{d_i}}{1-t}
\end{displaymath}

\noindent where $d_1, \dots , d_n$ are the degrees of a set of
fundamental homogeneous invariant polynomials of $W$ ($d_1, \dots ,
d_n$ are independent, up to reordering, of the polynomials choosen). By \cite[page 291]{5}, $d_1, \dots, d_n = m,2m, \dots , (n-1)m,dn$ when $W = G(m,p,n)$.\\

\noindent Hence, if the dual representation of $( \{ \underline{ \mu }
\}, \epsilon)$ is $( \{ \underline{ \lambda } \}, \eta)$, equation (\ref{eq:poin}) becomes

\begin{displaymath}
P(L(\{ \underline{\mu} \},\epsilon),t) =
\end{displaymath}

\begin{displaymath}
\frac{\dim ( \{ \underline{\mu} \},\epsilon) \, t^{b_{ \{ \underline{ \lambda } \} }} \, (1 - t^m)(1 - t^{2m}) \cdots (1 - t^{(n-1)m})(1 - t^{nd}) \, \prod_{i=0}^{m-1} H_{\lambda^i}(t^m) (1 - t^{mn})}{(1-t)^n(1 - t^{dn}) \, R_{ \{ \underline{\lambda} \}}(t)(t^m)_{(n)} \, \prod_{i = 0}^{m-1} t^{n(\lambda^i)m} } 
\end{displaymath}

\begin{equation}\label{eq:red2}
= \frac{\dim ( \{ \underline{\mu} \},\epsilon) \,  t^{b_{ \{ \underline{\lambda} \} }} \, \prod_{i=0}^{m-1} H_{\lambda^i}(t^m)}{(1 - t)^n \, R_{ \{ \underline{\lambda} \}}(t) \, \prod_{i = 0}^{m-1} t^{n(\lambda^i)m}}
\end{equation}

\noindent Let $k \in \N$ such that $t^k \mid R_{ \{
  \underline{\lambda} \}}(t)$ but $t^{k+1} \nmid R_{ \{
  \underline{\lambda} \}}(t)$ in $\Z[t]$ and write $R_{ \{ \underline{\lambda}
  \}}(t) = t^k \tilde{R}_{ \{ \underline{\lambda} \}}(t)$. Then
  rearrange equation (\ref{eq:fake}) as

\begin{equation}
f_{\{ \underline{\lambda} \} }(t) = \left(t^k \prod_{i = 0}^{m-1}
t^{n(\lambda^i)m}\right)\tilde{R}_{\{
  \underline{\lambda} \} }(t) \left( \frac{1 - t^{dn}}{1 - t^{mn}}(t^m)_{(n)} \prod_{i = 1}^m
\frac{1}{H_{\lambda^i} (t^m)}\right)
\end{equation}

\noindent Since each $H_{\lambda^i}(t^m)$ is a product of factors of the form $(1 - t^l)$, the product in the right most bracket consists entirely of factors of the form $(1 - t^l)$. Therefore  

\begin{displaymath}
t^{b_{ \{ \underline{\lambda} \} }} = t^k \prod_{i = 0}^{m-1} t^{n(\lambda^i)m}
\end{displaymath}

\noindent and equation (\ref{eq:red2}) becomes

\begin{equation}\label{eq:redbig}
P(L(\{ \underline{\mu} \},\epsilon),t) = \frac{\dim ( \{
  \underline{\mu} \},\epsilon) \prod_{i = 1}^m H_{\lambda^i}(t^m)}{(1
  - t)^n \tilde{R}_{ \{ \underline{\lambda} \} }(t)}.
\end{equation}

\noindent To contradict Lemma \ref{lem:poly} and hence prove Proposition \ref{prop:mainprop} we have 

\begin{lem}
Let $p \neq 1$ and $W = G(m,p,n)$ with $W \neq G(2,2,3)$. Then there exists $( \{ \underline{\mu} \},\epsilon) \in \Irr(W)$ such that the right hand side of equation (\ref{eq:redbig}) is not an element of $\C[t]$.
\end{lem}

\begin{proof}We consider the cases $n = 2,3$ and $n > 3$ separately. For $n > 3$ choose $( \{ \underline{\mu} \},\epsilon)$ such that its dual representation is $\underline{\lambda} = (\lambda^0, \emptyset, \dots \emptyset )$, where $\lambda^0 = (2,2,1,1, \dots 1)$. Then 

\begin{displaymath}
\tilde{R}(t) = R(t) = 1 + t^{dn} + t^{2dn} + \dots + t^{(p-1)dn} = \frac{1 - t^{mn}}{1 - t^{dn}}
\end{displaymath}

\noindent and for this particular $m$-multipartition we have

\begin{displaymath}
\prod_i H_{\lambda^i}(t^m) = H_{\lambda^0}(t^m) = (1 - t^{2m})(1 - t^m)(1 - t^{(n-1)m})(1 - t^{(n-2)m})\prod_{i = 1}^{n-4} (1-t^{im})
\end{displaymath}

\noindent Equation (\ref{eq:redbig}) becomes

\begin{equation}\label{eq:final}
P(L( \{ \underline{\mu} \},\epsilon),t) = \frac{\dim ( \{ \underline{\mu} \},\epsilon) (1 - t^{2m})(1 - t^m)(1 - t^{(n-1)m})(1 - t^{(n-2)m})(t^m)_{n-4} (1 - t^{dn})  }{(1 - t^{mn})(1-t)^n}.
\end{equation}

\noindent The numerator of (\ref{eq:final}) factorises over $\C$ as a product of factors $(1 - \omega t)$, where $\omega$ is a primitive $k^{th}$ root of unity with $1 \le k < mn$, whereas the denominator contains at least one factor of the form $(1 - \sigma t)$, where $\sigma$ is a primitive $(mn)^{th}$ root of unity. Therefore, since $\C[t]$ is an Euclidean domain, the right hand side of (\ref{eq:final}) cannot not lie in $\C [t]$.\\

\noindent For $n = 2$ and $m \ge n$, take $\underline{\lambda} = ((1), (1), \emptyset \dots \emptyset )$. Then

\begin{displaymath}
\prod_i H_{\lambda^i}(t^m) = (1 - t^m)^2 \qquad R(t) = \frac{t(1 - t^{2m})}{1 - t^{2d}}, \quad  \textrm{ and } \quad \tilde{R}(t) = \frac{1 - t^{2m}}{1 - t^{2d}}.
\end{displaymath}

\noindent Substituting into (\ref{eq:redbig}) 

\begin{displaymath}
P(L( \{ \underline{\mu} \},\epsilon),t) = \frac{\dim ( \{ \underline{\mu} \},\epsilon) (1 - t^m)^2 (1 - t^{2d})}{(1 - t^{2m})(1 - t)^2}.
\end{displaymath}

\noindent By the same reasoning as above, since $2m > 2d, m$, this rational function is not a polynomial.\\

\noindent Similarly, for $n = 3$ and $m \ge n$, take $\underline{\lambda} = ((1), (1), (1), \emptyset \dots \emptyset )$. Then

\begin{displaymath}
\prod_i H_{\lambda^i}(t^m) = (1 - t^m)^3 \qquad R(t) = \frac{t^3(1 - t^{3m})}{1 - t^{3d}}, \quad \textrm{ and } \quad  \tilde{R}(t) = \frac{1 - t^{3m}}{1 - t^{3d}}.
\end{displaymath}

\noindent Substituting into (\ref{eq:redbig}) 

\begin{displaymath}
P(L( \{ \underline{\mu} \},\epsilon),t) = \frac{\dim ( \{ \underline{\mu} \},\epsilon) (1 - t^m)^3 (1 - t^{3d})}{(1 - t^{3m})(1 - t)^3}.
\end{displaymath}

\noindent Once again, this rational function is not a polynomial because $3m > 3d, m$.
\end{proof}

\subsection{The Exceptional Groups}

Using the computer algebra program \cite{2} together with the package \cite{3} we calculate for each exceptional complex reflection group $W$ (excluding $G_4$), the number of irreducible representations $\lambda$ for which the polynomial $t^{-b_{\lambda^*}} f_{\lambda^*}(t)$ does not divide $P(\C[\frak{h}]^{coW},t)$ in $\C[t]$. Table (\ref{tab:tab1}) gives the results of these calculations. For each of these $\lambda$, Lemma \ref{lem:poly} does not hold and hence $\dim L(\lambda) < |W|$ for all values of $\bf{c}$. Since there is always at least one such $\lambda$ for every exceptional group, Proposition \ref{prop:mainprop} is proved for the exceptional groups.

\begin{table}[h]\label{tab:tab1}
\centering
\caption{Number of irreducibles that fail Lemma \ref{lem:poly}}
\begin{tabular}{c|ccccccccccccccccccc}
Group & 5 & 6 & 7 & 8 & 9 & 10 & 11 & 12 & 13 & 14 & 15 & 16 & 17 & 18 & 19 & 20 & 21 & 22 & 23 \\
\hline\hline
\# failures & 3 & 6 & 13 & 2 & 16 & 15 & 43 & 1 & 4 & 9 & 18 & 15 & 55 & 70 & 164 & 18 & 42 & 12 & 4 
\end{tabular}
\begin{tabular}{c|cccccccccccccc}
Group & 24 & 25 & 26 & 27 & 28 & 29 & 30 & 31 & 32 & 33 & 34 & 35 & 36 & 37 \\
\hline\hline
\# failures & 8 & 3 & 10 & 26 & 5 & 24 & 24 & 40 & 33 & 30 & 148 & 9 & 30 & 75
\end{tabular}
\end{table}

\noindent The code used to produce the data in Table (\ref{tab:tab1}) is available on the author's website \cite{10}. For every exceptional group, the fake polynomials of the irreducible characters are listed there. The remainder of $P(\C[\frak{h}]^{coW},t)$ on division by $t^{-b^*}f_{\lambda^*}(t)$ is also listed. In addition, this information is available for many of the groups $G(m,p,n)$ of rank $\le 5$. 

\section{The exceptional group $G_4$}

\noindent The group $G_4$, as labelled in \cite{5}, is the binary tetrahedral
group. It can be realised as a finite subgroup of the group of units in the quaternions

\begin{displaymath}
G_4 = \{ \pm 1,\pm i, \pm j , \pm k, \frac{1}{2}(\pm 1 \pm i \pm j
\pm k) \}
\end{displaymath}

\noindent and has order 24. It is generated by the elements $s_1 =
\frac{1}{2}(-1 + i + j - k)$ and $s_2 =  \frac{1}{2}(-1 + i - j + k)$
and has presentation $G_4 = \langle s_1,s_2 | s_1^3 = s_2^3 =
(s_1s_2)^6 = 1 \rangle$. It has seven conjugacy classes which we label $Cl_1 = \{
1 \}$, $Cl_2$, $Cl_3$, $Cl_4$, $Cl_5$, $Cl_6$, and $Cl_7$. The
character table is

\begin{displaymath}
\begin{array}{c|ccccccc}
\textrm{Class} & 1 & 2 & 3 & 4 & 5 & 6 & 7\\
\textrm{Size} & 1 & 1 & 4 & 4 & 6 & 4 & 4\\
\textrm{Order} & 1 & 1 & 3 & 3 & 4 & 6 & 6\\
\hline
T   & 1 & 1 & 1 & 1 & 1 & 1 & 1 \\
V_1 & 1 & 1 & \omega^2 & \omega & 1 & \omega^2 & \omega \\
V_2 & 1 & 1 & \omega & \omega^2 & 1 & \omega & \omega^2 \\
W  & 2 & -2 & -1 & -1 & 0 & 1 & 1\\
\frak{h} & 2 & -2 & -\omega^2 & -\omega & 0 & \omega^2 & \omega \\
\frak{h}^* & 2 & -2 & -\omega & -\omega^2 & 0 & \omega & \omega^2 \\ 
U & 3 & 3 & 0 & 0 & -1 & 0 & 0
\end{array}
\end{displaymath}

\noindent where $\omega$ is a primitive cube root of unity. Note that
the reflection representation $\frak{h}$ has dimension $2$, therefore
$G_4$ is a rank $2$ complex reflection group. \\

\noindent The group $G_4$ has two classes which consist of complex reflections and we label
these reflections as 

\begin{displaymath}
Cl_3 = \{s_1,s_2,s_3,s_4 \}
\end{displaymath}
\begin{displaymath}
= \{ \frac{1}{2}(-1 + i + j - k), \frac{1}{2}(-1 + i - j + k),
\frac{1}{2}(-1 - i + j + k), \frac{1}{2}(-1 - i - j - k) \}
\end{displaymath}
\noindent and 
\begin{displaymath}
Cl_4 = \{t_1,t_2,t_3,t_4 \}
\end{displaymath}
\begin{displaymath}
= \{\frac{1}{2}(-1 - i - j + k) , \frac{1}{2}(-1 + i - j - k),
\frac{1}{2}(-1 - i + j - k), \frac{1}{2}(-1 + i + j + k)\}
\end{displaymath}

\noindent Unlike all other exceptional irreducible complex reflection groups we have

\begin{thm}\label{thm:G4}
For generic values of $\mathbf{c}$, the generalised Calogero-Moser space $X_{\bf{c}}$ associated to $G_4$ is a smooth variety.
\end{thm}

\begin{proof}
The theorem is proved by showing that each irreducible $H_{0,\mathbf{c}}$-module is isomorphic to the regular representation of $G_4$. By \cite[Proposition 3.8]{1}, this is equivalent to the statement of the theorem. Let $E = T \oplus V_1 \oplus V_2 \oplus 3U$ and $F = \frak{h} \oplus \frak{h}^* \oplus W$ be two $G_4$-modules.\\

\subsection*{Claim 1} Let $L$ be a finite dimensional $H_{0,\mathbf{c}}$-module, then $L \cong aE \oplus bF$, for some $a,b \in \Z_{\ge 0}$.\\

\noindent To prove Claim 1 we use an argument similar to that of \cite[Proposition 16.5]{1}. Let $\rho : H_{0,\mathbf{c}} \rightarrow \End_{\C}(L)$ realise the action of $H_{0,\mathbf{c}}$ on $L$. Then, for all $x,y \in \frak{h} \oplus \frak{h}^*$, we have the commutation relation

\begin{equation}\label{eq:comm}
[\rho(x),\rho(y)] = c_1 \sum_{i = 1}^4 \omega_{s_i}(x,y) \, \rho(s_i) + c_2 \sum_{j = 1}^4
\omega_{t_j} (x,y) \, \rho(t_j) 
\end{equation}

\noindent By Lemma \ref{lem:nondeg}, $\sum_{i = 1}^4 \omega_{s_i} = \sum_{j = 1}^4 \omega_{t_j} = 2 \omega$. Taking traces on both sides of equation (\ref{eq:comm})

\begin{equation}\label{eq:trace}
0 = c_1 \, 2 \omega(x,y) \, Tr_L (s_1) + c_2 \, 2 \omega (x,y) \, Tr_L ( t_1)  \qquad \forall x,y \in \frak{h} \oplus \frak{h}^*\end{equation}

\noindent Since $c_1$ and $c_2$ are generic i.e. take values in a dense open subset of $\C^2$, and equation (\ref{eq:trace}) is linear, we have $0 = 2 \omega(x,y) \, Tr_L (s_1) = 2 \omega (x,y) \, Tr_L ( t_1)$. The fact that $\omega$ is
nondegenerate implies that $Tr_L$ is zero on $Cl_3$ and $Cl_4$.\\ 

Using the fact that $s_1$ is a complex reflection and $\dim
\frak{h}^* = 2$, we can choose a nonzero $x_1 \in \frak{h}^*$ such that
$s_1(x_1) = x_1$. Then $s_1[x_1 , y] = [x_1 , s_1y]$ for all $y \in
\frak{h}$.  Since $s_1(x_1) = x_1$, $x_1 \in \ker (1 - s_1)$ and hence $\omega_{s_1}(x_1,y) = 0$ for all $y \in \frak{h}$.  Similarly, $s_1t_1 = 1$ implies that $x_1 \in Fix(t_1)$ and hence $\omega_{t_1}(x_1,y) = 0$. Therefore, multiplying both sides of equation (\ref{eq:comm}) on the left by $\rho(s_1)$ and taking traces

\begin{displaymath}
0 = c_1 \sum_{i = 2}^4 \omega_{s_i}(x_1,y) \, Tr_L(s_1s_i) + c_2 \sum_{j = 2}^4 \omega_{t_j} (x_1,y) \, Tr_L(s_1t_j)
\end{displaymath}

\noindent Again, using the fact that $c_1, c_2$ are generic, we get  
\begin{displaymath}
0 = \sum_{i = 2}^4 \omega_{s_i}(x_1,y) \, Tr_L(s_1s_i) = \sum_{j = 2}^4 \omega_{t_j}(x_1,y) \, Tr_L(s_1t_j) 
\end{displaymath}

\noindent Since $s_1s_2, s_1s_3$ and $s_1s_4$ all belong to $Cl_7$ and $s_1t_2, s_1t_3,s_1t_4$ all belong to $Cl_5$ we have 
\begin{displaymath}
0 = \sum_{i = 2}^4 \omega_{s_i}(x_1,y)\, Tr_L(s_1s_i) = 2 \omega(x_1,y) \, Tr_L(s_1s_2)
\end{displaymath}
\begin{displaymath}
0 = \sum_{j = 2}^4 \omega_{t_j}(x_1,y) \, Tr_L(s_1t_j) = 2 \omega(x_1,y) \, Tr_L(s_1t_2)
\end{displaymath}
Therefore $Tr_L$ is zero on $Cl_7$ and $Cl_5$.\\
We can also multiplying both sides of equation (\ref{eq:comm}) on the left by $\rho(t_1)$ instead of $\rho(s_1)$. Noting that $t_1^2 \in Cl_3$, $t_1t_2,t_1t_3,t_1t_4 \in Cl_6$ and repeating the above argument shows that $Tr_L$ is also zero on $Cl_6$.\\

\noindent Therefore any element of $G_4$ that has nonzero trace on $L$ must belong to $Cl_1$ or $Cl_2$. Hence the character associated to $L$ must take values $(n,m,0,0,0,0,0)$, for some $n \in \Z_{\ge 0}, m \in \Z$, on the classes $Cl_1, Cl_2, \dots , Cl_7$. Taking inner products shows that

\begin{displaymath}
L \cong \frac{1}{|G_4|}(n + m) E \oplus \frac{2}{|G_4|}(n - m)F
\end{displaymath}

\noindent Setting $a = \frac{1}{|G_4|}(n + m)$ and $b = \frac{2}{|G_4|}(n - m)$ proves Claim 1. \\

\subsection*{Claim 2} Let $L$ be an irreducible representation of $H_{0,\mathbf{c}}$, with $\mathbf{c}$ generic. Then $L$ must be isomorphic to $E \oplus F$ or $\C G_4$ as a $G_4$-module.\\

\noindent If $L$ is irreducible then $\dim L \le 24$. Therefore Claim 1 implies that $L \cong E,2E,nF,1 \le n \le 4, E \oplus F$ or $\C G_4$. Assume that $L$ is isomorphic to $E$ as a $G_4$-module. The action of $\frak{h}^*$ on $L$ defines a $G_4$-equivariant linear map $\phi : \frak{h}^* \rightarrow
\End_{\C} (E)$. The $G_4$-module $\End_{\C}(E)$ decomposes as 

\begin{displaymath}
\End_{\C}(E) \cong (T \otimes T) \oplus 2(T \otimes V_1) \oplus  2(T \otimes V_2) \oplus 6(T \otimes U) \oplus (V_1 \otimes V_1) \oplus 2(V_1 \otimes V_2) \oplus
\end{displaymath}
\begin{displaymath}
6(V_1 \otimes U) \oplus (V_2 \otimes V_2) \oplus 6(V_2 \otimes U) \oplus 9(U \otimes U) \cong 12T \oplus 12V_1 \oplus 12V_2 \oplus 36U
\end{displaymath}

\noindent This shows that $\frak{h}^*$ is not a summand of $\End_{\C}(E)$. Thus $\phi$ must be the zero map. Similarly, the action of $\frak{h}$ must also be zero on $E$. This implies that the right hand side of equation (\ref{eq:comm}) must also act as zero on $E$. In particular, it must act as zero on $T \subset E$. This means that 

\begin{displaymath}
0 = c_1 \sum_{i = 1}^4 \omega_{s_i}(x,y) + c_2 \sum_{j = 1}^4 \omega_{t_j}(x,y) = 2(c_1 + c_2)\omega(x,y)
\end{displaymath}

\noindent This is a contradiction because $c_1,c_2$ are generic and $\omega$ is nondegenerate. Hence $L$ cannot be isomorphic to $E$. Repeating the above argument for $F$ we have 

\begin{displaymath}
\End_{\C}(F) \cong (\frak{h} \otimes \frak{h}) \oplus 2(\frak{h} \otimes \frak{h}^*) \oplus 2(\frak{h} \otimes W) \oplus
\end{displaymath}
\begin{displaymath}
(\frak{h}^* \otimes \frak{h}^*) \oplus 2(\frak{h}^* \otimes W) \oplus (W \otimes W) \cong 3T \oplus 3V_1 \oplus 3V_2 \oplus 9U
\end{displaymath}

\noindent Therefore $\frak{h}^*$ and $\frak{h}$ must act as zero on $F$. If we consider the right hand side of equation (\ref{eq:comm}), this time restricted to $W \subset F$ then we have

\begin{displaymath}
0 = c_1 \sum_{i = 1}^4 \omega_{s_i}(x,y) \, \rho|_{W}(s_i) + c_2 \sum_{j = 1}^4 \omega_{t_j}(x,y) \, \rho|_W(t_j)
\end{displaymath}

\noindent Taking the trace of this equation gives $0 = -2(c_1 + c_2)\omega(x,y)$, which is a contradiction because $c_1,c_2$ are generic and $\omega$ is nondegenerate. Therefore $L \not\cong F$. The same reasoning shows that $L$ cannot be isomorphic to $2E$ or $nF, 2 \le n \le 4$ either. This proves Claim 2.\\

\subsection*{Claim 3} Let $L$ be an irreducible $H_{0,\mathbf{c}}$-module, then $L$ cannot be isomorphic to $E \oplus F$ as a $G_4$-module.\\

\noindent By Lemma \ref{lem:Verma}, there exists a generalised Verma module $\Delta_{\mathbf{c}}(M;p,q)$ and a surjective homomorphism $\phi : \Delta_{\mathbf{c}}(M;p,q) \twoheadrightarrow L$. As a $G_4$-module we have 

\begin{displaymath}
\Delta_{\mathbf{c}}(M;p,q)  = H_{0,\mathbf{c}}(W) \otimes_{\C[\frak{h}]^W \otimes
  \C[\frak{h}^*] \rtimes W_q} M \cong \C G_4 \otimes \textrm{Ind}_{(G_4)_q}^{G_4} M \cong k \C G_4
  \end{displaymath}
  
\noindent where $(G_4)_q$ is the stabiliser of $q \in \frak{h}^*$ and $k = \lbrack G_4 : (G_4)_q \rbrack \dim M$. The generalised Verma module $\Delta_{\mathbf{c}}(M;p,q)$ has a finite composition series. Each factor of this series must have dimension $\le 24$. Therefore, by Claim 2, each factor is isomorphic to either $\C G_4$ or $E \oplus F$ as a $G_4$-module. Hence there exist $m,n \in \N$ such that $k \C G_4 \cong m \C G_4 \oplus n (E \oplus F)$ with $n \ge 1$. But then $n (E \oplus F) \cong (k - m) \C G_4$, which is a contradiction. This completes the proof of Claim 3 and the theorem.
\end{proof}

We can now apply a result of Namikawa \cite[Corollary 2.10]{9}, which we state for the convenience of the reader.

\begin{thm}[Namikawa]
Let $(X,\omega)$ be an affine symplectic variety equipped with a $\C^*$-action such that
\begin{itemize}

\item the weights of $\C^*$ on $X$ are all positive and there exists a unique fixed point $0 \in X$,

\item the symplectic form $\omega$ has positive weight $l > 0$.

\end{itemize}
Then the following are equivalent
\begin{enumerate}

\item $X$ has a crepant resolution

\item $X$ has a smoothing by a Poisson deformation.

\end{enumerate}
\end{thm}

\begin{cor}\label{cor:res}
Let $X$ be the symplectic singularity $\mathfrak{h} \times \mathfrak{h}^* / G_4$. There exists a symplectic resolution $\pi : Z \rightarrow X$ of $X$.
\end{cor}

\noindent First we recall some definitions from \cite[page 236]{18}, the reader should consult that article for details. A variety will mean a quasi-projective variety over $\C$. Let $X,Y$ be normal varieties with $K_X$ $\Q$-Cartier and $f : Y \rightarrow X$ a birational morphism. We can write
\begin{displaymath}
K_Y \equiv f^*(K_X) + A
\end{displaymath} 
If $E$ is a prime exceptional divisor on $Y$ then the \textit{discrepancy} of $E$ with respect to $X$ (denoted $a(E,X)$) is defined to be the coefficient of $E$ in $A$. If $f' : Y' \rightarrow X$ is another birational morphism and $E' \subset Y'$ the birational transform of $E$ on $Y'$ then $a(E,X) = a(E',X)$. Therefore $a(E,X)$ depends only on $E$ and not on $Y$. The variety $X$ is called \textit{canonical} if $a(E,X) \ge 0$ for all $E$.

\begin{proof}
The affine variety $X$ is four dimensional and normal. By \cite[Watanabe's Theorem]{15} $X$ has Gorenstein singularities and hence the canonical divisor $K_X$ is trivial (and hence Cartier). The affine variety $V = \mathfrak{h} \times \mathfrak{h}^*$ is smooth and therefore $V$ is canonical. Since $G_4$ is a finite group, the quotient map $\pi : V \rightarrow X$ is a finite dominant morphism and $\pi^* K_X = \pi^* \mathcal{O}_X = \mathcal{O}_V = K_V$. Therefore we can apply \cite[Proposition 3.16]{18} which says that $X$ is canonical. Therefore the pair $(X,\emptyset)$ is a Kawamata log terminal pair (as defined in \cite{16}) and we can apply \cite[Lemma 2.1]{16} to conclude that there exists an effective $\Q$-factorial terminal pair $(Y,B)$ together with a birational morphism $f : Y \rightarrow X$ such that 
\begin{displaymath}
K_Y + B \equiv f^* (K_X)
\end{displaymath}
However as noted above we can write $K_Y \equiv f^*(K_X) + A$ with $A = - B$. Since $X$ is canonical $a(E,X) \ge 0$ for all exceptional prime divisors $E$ on $Y$. Hence $A$ is an effective divisor. But $B$ is also effective therefore $A = B = 0$ and we deduce that $f : Y \rightarrow X$ is a crepant morphism. As noted in \cite[Section 4.14]{1}, $\{ X_{\bf{c}} \}_{\bf{c} \in \C^2}$ is a Poisson deformation of $X$. Therefore Theorem \ref{thm:G4} says that $X$ has a smoothing by a Poisson deformation. Now we can apply Namikawa's result \cite[Theorem 2.4]{9} and conclude that there exists a symplectic resolution $\pi : Z \rightarrow X$.
\end{proof}

\section{Acknowledgements}

The research described here was done at the University of Edinburgh with the financial support of the EPSRC. This material will form part of the author's PhD thesis for the University of Edinburgh. The author would like to thank Iain Gordon for suggesting this problem and for his help, encouragement and patience. He would also like to thank Maurizio Martino for help with the proof of Lemma \ref{lem:nondeg} and many useful conversations and Baohua Fu and Ivan Cheltsov for assistance with Corollary \ref{cor:res}. The author would like to thank Ulrich Thiel for pointing out a mistake in the published version of the article.

\end{document}